\documentclass[11pt]{amsart}
\usepackage{epsfig,latexsym,amsfonts,amssymb,amsmath,amscd}
\usepackage{amsthm}
\usepackage[mathscr]{eucal}
\usepackage{mathrsfs}
\usepackage{oldgerm,units,color}

\makeatletter
\def\section{\@startsection{section}{1}%
  \z@{1.1\linespacing\@plus\linespacing}{.8\linespacing}%
  {\normalfont\Large\scshape\centering}}
\makeatother

\theoremstyle{plain}

\newtheorem*{conj*}{Root Groups Conjecture}
\newcommand{\etype}[1]{\renewcommand{\labelenumi}{(#1{enumi})}}
\newtheorem*{thm1.2}{(1.2) Theorem}
\newtheorem*{thm1.3}{(1.3) Theorem}
\newtheorem*{thm1.4}{(1.4) Theorem}
\newtheorem*{prop*}{Proposition}

\def\diag{\operatorname{diag}}
\def\lran{\operatorname{lrank}}
\def\wk{{\operatorname{wk}}}
\def\rran{\operatorname{rrank}}
\def\rank{\operatorname{rank}}
\def\Miy{\operatorname{Cl}}
\def\adiag{\operatorname{antidiag}}
\def\fun{cross\ }
\def\comm{{\operatorname {comm} }}
\def\ba{\bar a}
\def\mcE{\mathcal E}
\def\tmcE{\tilde{\mathcal E}}
\def\mcI{\mathcal I}
\def\mcA{\mathcal A}
\def\mcM{\mathcal M}
\def\mcP{\mathcal P}
\def\mcS{\mathcal S}
\def\bb{\bar b}
\def\bc{\bar c}
\def\bd{\bar d}
\def\by{\bar y}
\def\Af{A^{K}}
\def\grA{{A_{\operatorname{gr}}}}
\def\pp{{\operatorname{preprim}}}
\def\npp{{\operatorname{nspreprim}}}
\newtheorem{prop}{Proposition}[section]

\newtheorem{ques}[prop]{Question}
\newtheorem{cor}[prop]{Corollary}
\newtheorem{lemma}[prop]{Lemma}

\theoremstyle{definition}
\newtheorem{Def}[prop]{Definition}
\newtheorem*{Def*}{Definition}
\newtheorem{Defs}[prop]{Definitions}

\newtheorem{Defsnot}[prop]{Definitions and notation}
\newtheorem{example}[prop]{Example}

\newtheorem*{notation*}{Notation}
\newtheorem{remark}[prop]{Remark}

\newtheorem*{remark*}{Remark}

\DeclareMathOperator{\tld}{\sim\!}
\DeclareMathOperator{\inv}{inv}

\newcommand{\mouf}{\mathbb{M}}
\newcommand{\mon}{\operatorname{mon}}
\newcommand{\cala}{\mathcal{A1}}
\newcommand{\calb}{\mathcal{B}}
\newcommand{\calc}{\mathcal{C}}
\newcommand{\cald}{\mathcal{D}}
\newcommand{\cale}{\mathcal{E}}
\newcommand{\calf}{\mathcal{F}}
\newcommand{\calg}{\mathcal{G}}
\newcommand{\calh}{\mathcal{H}}
\newcommand{\calj}{\mathcal{J}}
\newcommand{\calk}{\mathcal{K}}
\newcommand{\cali}{\mathcal{I}}
\newcommand{\call}{\mathcal{L}}
\newcommand{\calm}{\mathcal{M}}
\newcommand{\caln}{\mathcal{N}}
\newcommand{\calo}{\mathcal{O}}
\def\ba{\bar a}
\def\mcE{\mathcal E}
\def\mcB{\mathcal B}
\def\tmcE{\tilde{\mathcal E}}
\def\mcI{\mathcal I}
\def\mcA{\mathcal A}
\def\mcM{\mathcal M}
\def\mcN{\mathcal N}
\def\core{\operatorname{core}}
\def\mcP{\mathcal P}
\def\mcS{\mathcal S}
\def\mcT{\mathcal T}
\def\bb{\bar b}
\newcommand{\rad}{\operatorname{rad}}
\newcommand{\calot}{\tilde{\mathcal{O}}}
\newcommand{\calp}{\mathcal{P}}
\newcommand{\calq}{\mathcal{Q}}
\newcommand{\calr}{\mathcal{R}}
\newcommand{\cals}{\mathcal{S}}
\newcommand{\calt}{\mathcal{T}}
\newcommand{\calu}{\mathcal{U}}
\newcommand{\calv}{\mathcal{V}}
\newcommand{\calw}{\mathcal{W}}
\newcommand{\calx}{\mathcal{X}}
\newcommand{\calz}{\mathcal{Z}}
\newcommand{\solv}{\mathcal{SLV}}
\newcommand{\cc}{\mathbb{C}}
\newcommand{\ff}{\mathbb{F}}
\newcommand{\bbg}{\mathbb{G}}
\newcommand{\LL}{\mathbb{L}}
\newcommand{\kk}{\mathbb{K}}
\newcommand{\mm}{\mathbb{M}}
\newcommand{\bmm}{\bar{\mathbb{M}}}
\newcommand{\nn}{\mathbb{N}}
\newcommand{\oo}{\mathbb{O}}
\newcommand{\mbp}{\mathbb{P}}
\newcommand{\qq}{\mathbb{Q}}
\newcommand{\rr}{\mathbb{R}}
\newcommand{\z}{\mathbb{U}}
\newcommand{\vv}{\mathbb{V}}
\newcommand{\zz}{\mathbb{Z}}
\newcommand{\Cl}{\operatorname{Cl}}
\newcommand{\fraka}{\mathfrak{A1}}
\newcommand{\frakm}{\mathfrak{m}}
\newcommand{\frakB}{\mathfrak{B}}
\newcommand{\frakC}{\mathfrak{C}}
\newcommand{\frakh}{\mathfrak{h}}
\newcommand{\frakI}{\mathfrak{I}}
\newcommand{\frakM}{\mathfrak{M}}
\newcommand{\N}{\mathbb{N}}
\newcommand{\frakO}{\mathfrak{O}}
\newcommand{\frakP}{\mathfrak{P}}
\newcommand{\T}{\mathfrak{T}}
\newcommand{\frakU}{\mathfrak{U}}
\newcommand{\frakV}{\mathfrak{V}}

\newcommand{\ga}{\alpha}
\newcommand{\gb}{\beta}
\newcommand{\gc}{\gamma}
\newcommand{\gca}{\gamma_1}
\newcommand{\gcb}{\gamma_2}
\newcommand{\gcc}{\gamma_3}
\newcommand{\gC}{\Gamma}
\newcommand{\gCt}{\tilde{\Gamma}}
\newcommand{\gd}{\delta}
\newcommand{\one}{\mathbf 1}
\newcommand{\zero}{\mathbf 0}
\newcommand{\gre}{\epsilon}
\newcommand{\gve}{\varepsilon}
\newcommand{\gl}{\lambda}

\newcommand{\gn}{\nu}
\newcommand{\gro}{\omega}

\newcommand{\gvp}{\varphi}
\newcommand{\gr}{\rho}
\newcommand{\gs}{\sigma}
\newcommand{\gS}{\Sigma}
\newcommand{\gt}{\tau}
\newcommand{\gta}{{\tau_a}}
\newcommand{\gtb}{{\tau_b}}
\newcommand{\gtc}{{\tau_c}}
\newcommand{\gth}{\theta}
\newcommand{\gTH}{\Theta}

\newcommand{\rmA}{\mathrm{A1}}
\newcommand{\rmD}{\mathrm{D}}
\newcommand{\rmI}{\mathrm{I}}
\newcommand{\rmi}{\mathrm{i}}
\newcommand{\rmL}{\mathrm{L}}
\newcommand{\rmO}{\mathrm{O}}
\newcommand{\rmS}{\mathrm{S}}
\newcommand{\rmU}{\mathrm{U}}
\newcommand{\rmV}{\mathrm{V}}

\newcommand{\nsg}{\trianglelefteq}
\newcommand{\rnsg}{\trianglerighteq}
\newcommand{\At}{A^{\times}}
\newcommand{\Dt}{D^{\times}}
\newcommand{\Ft}{F^{\times}}
\newcommand{\Kt}{K^{\times}}
\newcommand{\kt}{k^{\times}}
\newcommand{\ktv}{k_v^{\times}}
\newcommand{\Pt}{P^{\times}}
\newcommand{\Qt}{Q^{\times}}
\newcommand{\calqt}{\calq^{\times}}

\newcommand{\Rt}{R^{\times}}
\newcommand{\Lt}{L^{\times}}

\newcommand{\charc}{{\rm char}}
\newcommand{\df}{\stackrel{\text{def}} {=}}
\newcommand{\lr}{\longrightarrow}
\newcommand{\ra}{\rightarrow}
\newcommand{\hr}{\hookrightarrow}
\newcommand{\lra}{\longrightarrow}
\newcommand{\sminus}{\smallsetminus}
\newcommand{\lan}{\langle}
\newcommand{\ran}{\rangle}
\newcommand{\Ab}{{\rm Ab}}
\newcommand{\Aut}{{\rm Aut}}
\newcommand{\End}{{\rm End}}
\newcommand{\Ker}{{\rm Ker\,}}
\newcommand{\im}{{\rm Im\,}}
\newcommand{\Hom}{{\rm Hom}}
\newcommand{\FRAC}{\rm FRAC}
\newcommand{\Nrd}{\rm Nrd}
\newcommand{\PSL}{\rm PSL}
\newcommand{\SL}{{\rm SL}}
\newcommand{\GL}{{\rm GL}}
\def\eroman{\etype{\roman}}
\def\atb{a^{\tau_b}}
\def\hab{\widehat \atb}
\def\hac{\widehat {a^{\tau_c}}}
\def\supp{\operatorname{supp}}
\def\hda{\widehat {d^{\tau_a}}}
\def\hdb{\widehat {d^{\tau_b}}}
\def\habc{\widehat {a^{\tau_b \tau_c}}}
\def\hbc{\widehat {b^{\tau_c}}}
\def\bta{b^{\tau_a}}
\def\hba{\widehat \bta}
\def\ha{\hat a}
\def\hX{\hat X}
\def\hb{\hat b}
\def\hc{\hat c}
\def\hd{\hat d}
\def\hx{\hat x}
\def\Eref{Example~\ref}
\def\Tref{Theorem~\ref}
\newcommand{\GF}{{\rm GF}}
\newcommand{\R}{{\rm R}}
\newcommand{\Sz}{\rm Sz}
\newcommand{\St}{\rm St}
\newcommand{\Sym}{{\rm Sym}}
\newcommand{\tr}{{\rm tr}}
\newcommand{\tor}{{\rm tor}}
\newcommand{\PGL}{{\rm PGL}}

\newcommand{\ch}{\check}
\newcommand{\s}{\star}
\newcommand{\bu}{\bullet}
\newcommand{\dS}{\dot{S}}
\def\lran{\operatorname{lrank}}
\def\wk{{\operatorname{wk}}}
\def\rran{\operatorname{rrank}}
\def\rank{\operatorname{rank}}
\def\Miy{\operatorname{Cl}}
\def\adiag{\operatorname{antidiag}}
\def\fun{cross\ }
\def\comm{{\operatorname {comm} }}
\def\ba{\bar a}
\def\mcE{\mathcal E}
\def\tmcE{\tilde{\mathcal E}}
\def\mcI{\mathcal I}
\def\mcA{\mathcal A}
\def\mcM{\mathcal M}
\def\mcP{\mathcal P}
\def\mcS{\mathcal S}
\def\bb{\bar b}
\def\bc{\bar c}
\def\bd{\bar d}
\def\by{\bar y}
\def\Af{A^{K}}
\def\grA{{A_{\operatorname{gr}}}}
\def\Rad{{\operatorname{Rad}}}
\def\pp{{\operatorname{preprim}}}
\def\npp{{\operatorname{nspreprim}}}\newcommand{\onto}{\twoheadrightarrow}
\newcommand{\HH}{\widebar{H}}
\newcommand{\NN}{\widebar{N}}
\newcommand{\GG}{\widebar{G}}
\newcommand{\gCgC}{\widebar{\gC}}

\newcommand{\hG}{\widehat{G}}
\newcommand{\hN}{\widehat{N}}
\newcommand{\hgC}{\widehat{\gC}}
\newcommand{\linv}{\varprojlim}
\newcommand{\In}{{\rm In}_K}
\newcommand{\Inc}{{\rm Inc}_K}
\newcommand{\Inv}{{\rm Inv}}
\newcommand{\tit}{\textit}
\newcommand{\tbf}{\textbf}
\newcommand{\tsc}{\textsc}
\newcommand{\hal}{\frac{1}{2}}
\newcommand{\half}{\frac{1}{2}}
\newcommand{\restr}{\upharpoonright}
\newcommand{\rmk}[1]{\noindent\tbf{#1}}
\newcommand{\widebar}[1]{\overset{\mskip1mu\hrulefill\mskip1mu}{#1}
                \vphantom{#1}}
\newcommand{\widedots}[1]{\overset{\mskip1mu\dotfill\mskip1mu}{#1}
                \vphantom{#1}}
\newcommand{\llr}{\Longleftrightarrow}

\newcommand{\tem}{{\bf S}}
\newcommand{\tnem}{{\bf NS}}

\numberwithin{equation}{section}

\begin{document}
\title[Axial radical]{The axial radical}
\author[Louis Halle Rowen ]
{Louis Halle Rowen$^*$}

\address{Louis Rowen\\
       Department of Mathematics\\
        Bar-Ilan University\\
        Ramat Gan\\
        Israel}
\email{rowen@math.biu.ac.il}
\thanks{The author was supported by the Anshel Pfeffer Chair}





\begin{abstract}
An {\it axial algebra}   is an algebra generated by a set $X$ of semisimple idempotents satisfying a prescribed set of fusion rules.  Some axial algebras (such as primitive axial algebras of Jordan type $\half$) have a Frobenius form, in which case one can define the radical $\Rad(A)$ of the Frobenius form. Recently, Mamontov, S.~Shpectorov, and V.~Zhelyabin  defined a radical $R(A;X)$ as the largest ideal not containing any axis of $X$, which is always contained in $\Rad(A)$, and proved that $R(A;X) = \Rad(A)$ when the Frobenius form is anisotropic.  In this case they  also proved that $A/R(A;X)$ is semisimple, yielding an elegant structure theory for primitive axial algebras with an anisotropic Frobenius form. Our objective is to recover these results, removing the assumption of a Frobenius form.  We provide alternative characterizations of the radical, including  a Jacobson-type radical (not the same as their Jacobson radical). We also replace the Frobenius form by a more general notion of orthogonality. Finally, we introduce a Baer-type radical and compare it to $R(A;X)$ and $\Rad(A)$.
\end{abstract}
\date{\today}

 \maketitle


\section{Introduction}\label{sec1}

Axial algebras --- algebras generated by semisimple idempotents --- were studied long ago by Peirce (in the associative case) and by Albert \cite{A2,A3,A4}. Axial algebras have become an active area of research  since the discovery of their relation with vertex operator algebras, 3-transposition groups and the classification of finite simple nonabelian groups,  and vector flows on manifolds, as described in the excellent overview by McInroy~\cite{M}.

Recently, attention has focused on describing the structure of axial algebras in terms of radicals. One particularly successful approach in the case of the existence of an anisotropic Frobenius form is given in  \cite{MSZ}. The radical theory of Jordan algebras is an important special case which has been studied extensively for several decades, cf.~\cite{J2,McC,HM}.

This paper, inspired by the superb set of questions in \cite{GS}, is motivated by the desire to understand the structure of axial algebras which do not necessarily have an anisotropic Frobenius form, in terms of various notions of a radical, recovering and extending \cite{MSZ}.

\section{Fundamental notions}


In this paper, $C$ always denotes a commutative Noetherian integral domain  containing an element~$\half$. (In the literature, $C$ would be
 an infinite field~$F$, but the more general situation is useful for handling generic algebras.) We assume throughout that we are given a finite set $S\subset C$ containing $0,1,$ such that $\mu -\nu$ are invertible in $C$ for each distinct $\mu, \nu\in S.$ In particular, when $S = \{0,1,\gr\},$ we assume that $\gr (1-\gr)$ is invertible in $C$. 
$A$ always is assumed to be a faithful commutative (not necessarily associative) $C$-algebra.

\begin{Defsnot}\label{not1}
 We review some familiar definitions.
\begin{enumerate}

\item
For $a\in A,$ write $L_a$ for the left multiplication map
$z\mapsto az$.  We write $L_a^k(y)$ for $a(a(\dots(ay))),$ taken $k$ times.

\item For $\gr\in S,$ the $\gr$-\textbf{eigenspace} $A_\gr (a) $ of~$a$ is $ \{y\in A: L_ay = \gr y\}.$ (We permit the possibility that $A_\gr(a)=0.$)
 \item
 An \textbf{axis}
$a$ (with respect to $S$) is a nonzero
idempotent for which $L_a$ is semisimple, by which we mean  $\prod _{\gl \in S}(L_a-\gl)=0.$ ($\gl$ means $\gl$ times the identity map.)

 Obviously $1\in S$ since $L_a(a)=a^2 =a.$

\item  The idempotent $a$ is \textbf{primitive} if $A_1(a) = C a;$ in this case we write $y_1 = \gvp_a(y) a$ for $\gvp_a(y)\in C.$

 \item An algebra $A$ is \textbf{axial over} a countable\footnote{One could take $X$ arbitrary, but then would occasionally encounter issues of transfinite induction.} set $X$ if $A$ is generated by a set $X$ of axes.  Sometimes, we suppress $X$ in the notation, but many results depend on the choice of $X$. 
 in this case we write $X = \{ a_i : 1\le i \le N\}.$

 \item An axial algebra $A$   over $X$ is \textbf{primitive} if each axis in $X$ is primitive.

 \item The notation $(a_1\dots a_t)y$ for $t>2$ means $L_{a_1}L_{a_2}\dots L_{a_t}(y)$ i.e., inductively, $a_1 ((a_2 \dots a_t)y).$   $a_1\dots a_t$ denotes $(a_1\dots a_{t-1})a_t.$

 \item A \textbf{monomial} is a product $a_1\dots a_t$ for $a_i\in X.$
\end{enumerate} \end{Defsnot}

  \begin{remark}
    Every element of $A$ is a linear combination of products of monomials, so $A$ is countably spanned.
 \end{remark}

\begin{remark}  \textbf{Fusion rules} are rules of the form $A_\mu (a) A_{\nu} (a) \subseteq  \sum_{\gr \in S_{\mu,\nu}} A_\gr(a)$ for all $a\in X,$ with  $S_{\mu,\nu} \subseteq S$.

 The \textbf{basic fusion rules} for an axis $a$ are the three rules, writing $A_{0,1} = A_0+A_1$:
\begin{enumerate}  \item $A_{0}(a)^2 \subseteq A_0(a);$
 \item $A_{0,1}(a)$ is a $C$-subalgebra of $A$;
\item
 $A_{0,1}(a)A_{\gr} (a) \subseteq A_{\gr} (a)$  for each $\gr\in S$. In other words, $S_{0,\gr} = S_{1,\gr}=
 \gr$ for all $\gr\ne 0,1.$
\end{enumerate}
\end{remark}

\begin{example} Here are some well-known fusion rules. The tables are  accessible in~\cite{GS}.
    \begin{enumerate}\eroman
        \item An associative axial algebra satisfies the \textbf{Peirce fusion rules}, which are the basic fusion rules for $S = \{0,1\}.$

\item An \textbf{axial algebra of Jordan type $S$}, called an \textbf{axial algebra of Jordan type~$\gr$} when $S = \{0,1,\gr\},$ is an $S$-axial algebra for which each $a\in X $  satisfies $A_\mu(a) A_\nu(a) = 0$ for $\mu \ne \nu$ in $S\setminus \{1\}$, as well as  the basic fusion rules and the fusion rules $S_{\mu,\mu} \subseteq \{0,1\}$ for all $\mu\in S$.

 \item \cite{A1} All Jordan algebras generated by primitive idempotents are primitive of Jordan type $\half.$

\item An \textbf{axial algebra of Monster type} $(\ga,\gb)$ is a filtered algebra for $S = \{0,1,\ga,\gb\}.$
    \end{enumerate}
\end{example}

 \begin{remark}
 \cite{Row26} showed that the basic fusion rule (a) follows from (b) and (c) for every primitive axial algebra of Jordan type   over a field. We wonder if this holds for other primitive axial algebras.
 \end{remark}

\subsubsection{Applications of commutative ring theory}

For later reference, we recall some basic facts from commutative ring theory which impact on axial algebras.

\begin{remark}\label{fg} 
Let $\gl$ denote a commutative and associative indeterminate over $C.$
\begin{enumerate}\eroman
\item $C[\gl]$ is a Noetherian integral domain.

 \item For any axis $a\in A,$ $A$ can be viewed as a module over $C[\gl]$ via the action $f(\gl)y = f(L_a)(y)$. Thus, the theory of modules over commutative Noetherian rings can be applied.

\item When $C$ is a field, $C[\gl]$ is a principal ideal domain, over which any finitely generated module is a finite direct sum of cyclic modules, cf.~\protect{\cite[Theorem~2.54]{Row2}}.
 \end{enumerate}
\end{remark}

Remark~\ref{fg}(iii) simplifies the theory considerably for finite dimensional axial algebras over a field, but we can recover many results for arbitrary axial algebras.

\begin{prop} (Easy version of Chinese remainder theorem) Viewing $\gr$ as $\gr$ times the identity map, suppose $\prod _{\gr\in S}(L_a -\gr) = 0.$
    \begin{enumerate}\eroman
\item
$A =\bigoplus
_{\gr\in S}  A_{\gr} (a)$. Any element $y\in A$ is written uniquely as
\begin{equation}
 \label{dec2}
y = y_1+ \sum _{\gr \in S\setminus\{1\}} y_{\gr} , \qquad  \quad ay_\gr= \gr y_\gr. 
\end{equation}
\item $A_{\nu}(a) =  \prod _{\gr\ne \nu}(L_a -\gr) A,$ for each $\nu \in S.$
    \end{enumerate}
\end{prop}
\begin{proof} For any given $\nu \in S$, let $\gb_\nu = \prod _{\gr\ne \nu} ( \nu -\gr),$ which is invertible by hypothesis. Then $$\pi_\nu := \frac{1}{\gb_\nu}\prod _{\gr\ne \nu}(L_a -\gr) :A \to A_\nu$$
is a projection. Indeed, for each $y\in A,$ $(L_a -\nu)(\pi_\nu(y)) = \frac{1}{\gb_\nu}\prod _{\gr\in S}(L_a -\gr)(y) =0$, implying $\pi_\nu(y) \in A_\nu,$ whereas for $y\in A_\nu,$  $\gb_\nu = \prod _{\gr\ne \nu} ( \nu -\gr),$ $\pi_\nu(y) = \frac{\gb_\nu}{\gb_\nu} y = y.$ Hence we get (i), with the decomposition, so the projections $\pi_\nu$ are onto, yielding (ii).
\end{proof}

\begin{cor}
    If $\mcI \trianglelefteq A,$ then the image of $A_\gl$ in $\bar A =A/\mcI$ is $\bar A_\gl$.
 \end{cor}
\begin{proof}
    $\subseteq$ is immediate, but the sum is direct  since  taking components commutes with homomorphic images, by (ii).
\end{proof}

\begin{Def} An axial algebra $A$ satisfying the basic fusion rules is
    \begin{enumerate}
       \item  {\it graded} if $|S_{\mu,\nu}| =1$ for $\mu\ne \nu$ in $S$ and $S_{\mu,\mu} \subseteq \{0,1\} $ for each $\mu \in S$; in particular, $A$ is {\it 2-graded} with  $A^+= A_{0,1}(a)$  and $ A^- = \sum _{\gr \ne 0,1} A_\gr$. 

\item {\it filtered} if we can list the elements $S = \{0,1,\gr_1,\dots,\gr_t\}$, such that for each $a\in X,$ and   $S_{\gr_i,\gr_j}   =\{\gr_j \}$ whenever $i$ precedes $j$, and $S_{\gr_i,\gr_i} \subseteq \{\gr: \gr \text{ precedes } \gr_i\} .$
    \end{enumerate}
\end{Def}

\subsubsection{Miyamoto maps}

Here is one way of obtaining new axes.

\begin{Defs}\label{Miy7}
 \begin{enumerate}
 \item Suppose that $a$ is an $S$-axis of $A$. Picking nonempty $S' \subseteq S \setminus \{0,1\}$, define the map $\gt_{a,S'} :A \to A$ by
 $y \mapsto  \sum _{\gr \in S\setminus S'} {y_{\gr}}\, - \sum _{\gr \in S'} y_\gr$. $\gt_{a,S'},$ called a \textbf{Miyamoto map} of $a$,  clearly is a $C$-linear map of $A$, whose square is the identity map.

\item A \textbf{Miyamoto involution} is a Miyamoto map which is an
 algebra automorphism.
 \item A \textbf{Miyamoto axis} is an axis $a$ whose Miyamoto maps all are Miyamoto involutions. 
\end{enumerate}
\end{Defs}

As in \cite[Definition~4.1]{HSS2}, we note:

\begin{lemma} Any axial algebra $A$ of Jordan type is spanned by the images of axes of $X$ under Miyamoto involutions.
\end{lemma}
\begin{proof}(sketch) An easy observation on 2-graded algebras, noting that
 for any $S' \subseteq S\setminus \{0,1\}$, $A$ is 2-graded, where $A_{\gr}(a)$ is odd for $\gr \in S'$ and all other eigenspaces of $a$ are even; this shows that $ab = \sum \gl b_\gl$ is spanned by $a$ and Miyamoto images of $b$, and then one concludes by induction on the length of a word in the axes.
\end{proof}

The theory of algebras spanned by axes is particularly robust.

 \subsubsection{Frobenius forms}
 \begin{Defs}
 \begin{enumerate}
 \item An \textbf{associative bilinear form} $(\phantom{w},\phantom{w}): A \times A \to C$ satisfies
 $$\left(\sum \ga_{i }x_i,\sum \gb_{ j}y_j\right) = \sum _{i,j} \ga_i \gb_j (x_i,y_j), \quad (xy,z) = (x,yz), \qquad \forall x,y,z\in A.$$

 \item A \textbf{Frobenius form} is a symmetric associative bilinear form.\footnote{Any associative bilinear form on a commutative axial algebra is symmetric, cf.~\cite[Proposition~3.5]{HRS1}.}

\item An axis $a$ is \textbf{isotropic} if $(a,a)=0;$ otherwise $a$ is \textbf{anisotropic}. A Frobenius form is $X$-\textbf{anisotropic} if every  axis from $X$ is anisotropic. A Frobenius form is \textbf{normal} if $(a,a)=1$ for every primitive axis $a.$

\item The \textbf{radical} $\rad A$ of a Frobenius form on $A$, called the \textbf{Frobenius radical}, is $\{x\in A: (x,A)=0\}.$ 
 \end{enumerate}
 \end{Defs}

\begin{lemma}[{\cite[Lemma~2.7(ii)]{Row26}}] \label{normal} Suppose $a$ is an axis of an $S$-axial algebra $A$ with a Frobenius form, with $y,y' \in A$. Notation as in \eqref{dec2},
 \begin{enumerate}\eroman
 \item $(y_\mu,y'_{\nu})=0$ for all $\mu\ne \nu\in S$.

 \item If $a$ is primitive then $(a,a)y_1 = (a,y)a$. Thus, for $a$ anisotropic, and $C$ a field, $y_1 = \frac{(a,y)}{(a,a)}a$.

 \item  If $A$ is {\it spanned} by a set $X$ of primitive axes, then
 the following are equivalent:
 \begin{enumerate}
 \item $y \in \rad A$.

 \item $(a,y) =0$ for every $a\in X.$

 \end{enumerate}
\end{enumerate}
\end{lemma}

By \cite[Theorem~4.1]{HSS1}, every primitive axial algebra $A$ of Jordan type~$\gr$ over a field has a normal Frobenius form, given by $(a,y) = \gvp_a(y)$ for every primitive axis $a$ and all $y\in A$.

This also holds more generally for many primitive axial algebras of Jordan type~$S$, by \cite[Theorem~2.3]{RoSe6}, but a necessary condition for symmetry of the form clearly is $\gvp_a(b)= \gvp_b(a)$ for all primitive axes $a,b.$ Furthermore,
 the proof of \cite[Lemma~2.5]{RoSe6} utilizes Miyamoto involutions to prove the orthogonality of different components, a necessary condition by Lemma~\ref{normal}(i).

\begin{remark}\label{rady}
\begin{enumerate}\eroman
\item \cite{HRS1} studies those axial algebras of ``Monster type'' that have a Frobenius
form. As noted in \cite[Conjecture 8.15]{M}, Franchi, Mainardis, and Shpectorov conjecture that every axial algebra of Monster type has a Frobenius
form.
 \item A Frobenius form on $A$ passes to subalgebras,
 but need not pass to a homomorphic image $A/\mcI$ since it might not be well-defined. The Frobenius form is well defined on $A/\rad A.$
\end{enumerate}
 \end{remark}

 \subsubsection{Jordan algebras}

  Jordan algebras are defined by the identity
  $((xx)y)x = (xx)(yx)$. Axial Jordan algebras have Jordan type $\half$, in which all minimal nonzero idempotents are primitive, and have the bilinear operator
  $U_{a,b} = L_aL_b+L_bL_a -L_{ab}$ which yields the corresponding quadratic operator
  $U_a= U_{a,a}= 2L_a^2 - L_{a^2}.$
  Also, there is the Jordan triple product $\{a_1 a_2 a_3\}=(a_1a_2)a_3 +(a_2 a_3)a_1 -(a_1a_3)a_2$, which we  define in any axial algebra.

  \begin{remark}\label{Ua}
 For any  axis $a$ in an axial algebra of Jordan type $\half$ we would have $U_a(y) = 2a(ay)-ay = 2y_1 +\half y _\half -y_1 -\half  y_\half = y_1.$ In particular, if $a$ is primitive, $U_a(A) = Fa.$
  \end{remark}




 \section{$C[\gl]$-submodules of  $A$}

We start with the $C[\gl]$-module structure on $A,$ which depends on choosing $a\in X$.

 \begin{remark}\label{comp1} We view $\mcM_a:=A$ as a $C[\gl]$-module via the action $f(\gl)y = f(L_a)(y).$
     A key result is essentially \cite[Lemma~2.6]{MSZ}, which in our terminology says that   every  component of an element of a  $C[\gl]$-submodule $\mcN$ of $A$ is in $\mcN.$

     Accordingly we say a $C$-submodule $\mcN$ of $A$ is \textbf{$a$-homogeneous} if $y\in \mcN$ implies $y_\gl\in \mcN$ for all $\gl \in S.$

     It follows that if $y\in\mcN$ then $\tau(y)\in \mcN$ for any Miyamoto involution $\gt$ of $a$.
 \end{remark}

 \begin{Def}$ $
     \begin{enumerate}\eroman
         \item  Define the projection $\pi_a :A \to A_0(a)$ onto the 0 component by $\pi_a(y) = y_0.$

        \item  Define $A_\gr'(a) = \sum _{\mu \ne \gr} A_\mu(a),$ the complementary $C[\gl]$-submodule of  $A_\gr(a).$

     \end{enumerate}
 \end{Def}

\subsubsection{Modular submodules}

Modular left ideals are defined in \cite{J1} for associative algebras, as a way of avoiding a unit element. Quadratic ideals are defined in \cite{J2} for Jordan algebras, and play a key role in their structure theory.
Analogously, we say an $a$-homogeneous submodule $\mcN$ of $A$ is $a$-\textbf{modular} for $a\notin \mcN,$ if $ay-y\in \mcN$ for each $y\in \mcA.$
$\mcN$ is \textbf{$X$-axially modular} if there is an axis $a\notin \mcN$ in $X$ for which $\mcN$ is $a$-modular.

\begin{lemma}
    In an axial algebra $A$ over a field,   for any primitive  axis $a,$ the only $a$-modular submodule $\mcN$ of $A$ is  $ A_1'(a).$
\end{lemma}
\begin{proof} By hypothesis, $(L_a -1)(y) = ay-y\in \mcN$ for any $y\in \mcA. $ But $ (L_a -1)A =A_1'(a),$ implying $\mcN \supseteq A_1'(a)$. Furthermore
$A_1'(a)$ is $a$-modular since $ay - y = \sum _{\gr\in S} (\gr -1) y_\gr $ satisfies $\varphi_a(ay-y) = 0.$

  $A_1'(a)$ is maximal $a$-modular, since
if $z \in A \setminus A_1'(a)$ then $z_1 \ne 0,$ implying $a\in Ca$ and $Ca+A_1'(a)= A.$ Hence $\mcN = A_1'(a).$
\end{proof}

\subsection{Orthogonality}

Let us weaken the notion of Frobenius form, in order to see which hypotheses are needed in \cite[\S 5]{MSZ}.

\begin{Def}\label{orth}
    An \textbf{orthogonality relation}\footnote{We define a symmetric orthogonality relation; one could bypass symmetry, but the proofs would become more intricate.} $\perp$ on an axial algebra $A$ satisfies the following, writing   $a^\perp$ for $\{y\in A: a \perp y\},$  $V^\perp$ for $\{ w \in A: v\perp w\}$, and $V \perp W$ if $v\perp w$ for all $v\in V, w\in W$:
    \begin{enumerate}
    \item $A_\mu(a)\perp A_{\nu}(a)$ for  each $a\in X$, whenever $\mu \ne \nu.$ (Thus the orthogonality relation is nontrivial.)

\item $v\perp w$ if and only if $w\perp v$.


             \item If $ v_i \perp w$ for $v_i, v_2, w\in A$, then $\sum \ga_i v_i  \perp w$ for $\ga_i\in C.$

          \item If $u\perp yv$ for $y\in A,$ then $yu \perp v.$
    \end{enumerate}
\end{Def}

Heuristically, (1) provides some homogeneity, and (2),(3), and (4) say that orthogonality is commutative, bilinear, and associative. (This is reminiscent of the tensor product; the orthogonality relation can be viewed as a homomorphic image of $A\otimes_C A$.)

\begin{lemma}\label{perp3}
  \begin{enumerate}\eroman
     \item If each product of monomials is orthogonal to $y\in A,$ then $y\in A^\perp$.

     \item   $A^\perp \trianglelefteq A.$
  \end{enumerate}
\end{lemma}
\begin{proof} (i) By Definition~\ref{orth}(3),(4).

  (ii) $A^\perp$ is a $C$-submodule, by  Definition~\ref{orth}(2),(3). If $y\in A^\perp$ and $z\in A,$ then $Az \perp y,$ so $yz\perp A.$
\end{proof}

By definition, each Frobenius form defines an orthogonality relation, but  we want a canonical one (with respect to a generating set $X$). We shall always use  the minimal orthogonality relation, which could be obtained by taking the intersection of all orthogonality relations; it is obtained
 inductively, as follows:

\begin{Def}
At each stage, if $y\perp_i z$ we declare $y\perp_{i+1} z$, as well as the following:
  Define $\perp_1$ to satisfy  $A_\mu(a)\perp A_{\nu}(a)$ for  each $a\in X$, whenever $\mu \ne \nu.$ Given $\perp_{3i-2}$ for $i\ge 1,$ define $\perp_{3i-1}$ to satisfy   linearity, i.e., if $v \perp_{3i-2} w_i$ for $w_i\in A$, then $v \perp_{3i-1} \sum _i \ga_i w_i$  for $\ga_i\in C.$ Given $\perp_{3i-1}$,  define   $v\perp _{3i} yw$ if $yv \perp _{3i-1} w.$  Given $\perp_{3i}$, define   $v\perp _{3i+1} w$ if $w\perp _{3i} v$.

  Now declare $y\perp z$ when $y\perp_i z$ for some $i.$
\end{Def}

Obviously (1),(2), (3),   and   (4) hold by construction, and this is the minimal such orthogonality relation. Orthogonality has some advantages over a Frobenius form.

\begin{lemma}\label{perp0}
    Suppose   $\mcI \trianglelefteq A$.
    \begin{enumerate}\eroman
        \item   If $v\perp w$ in $A,$ then $\bar v\perp \bar w$ in $\bar A:= A/\mcI.$

            \item Conversely, for $\mcI \subseteq A^\perp,$ if  $\bar v\perp \bar w$ then $v\perp w$.
    \end{enumerate}
\end{lemma}
\begin{proof}
 (i)  The steps for building $\perp$ on $A$ are subsumed in the steps for building $\perp$ on $\bar A.$

 (ii)  The steps for building $\perp$ on $A$ are the same as the steps for building $\perp$ on $\bar A.$
 More explicitly, if $\bar v\perp \bar w$ then lifting the steps producing orthogonality shows $v'\perp w'$ for suitable respective representatives $v',w'$ of $\bar v, \bar w.$  Let  $y =v-v'\in A^\perp$ and $z =w-w'\in A^\perp$ yields $v = (y+v') \perp w'$ and certainly $v\perp z,$ and thus $v \perp (w'+z)=w.$
\end{proof}

\begin{lemma}\label{perp1}
    For $C$ a field, suppose $N$ is a $C[\gl]$-submodule of $A$. Then $a \perp N$ for any primitive axis $a\notin N$.
\end{lemma}
\begin{proof}
    For $y\in N,$ $a\perp y_\gl $ for all $\gl \ne 1,$ and $y_1 = 0$ since otherwise $a = \varphi_a^{-1}y_1 \in N.$ Hence $a \perp y.$
\end{proof}

\begin{Def}
    An axis $a$ is $\perp$-$\textbf{isotropic}$ if $a\perp A_1(a)$.
\end{Def}

(When $a$ is a primitive axis and $C$ is a field, this means $a\perp a$.)

\begin{lemma}\label{sing}
    An  axis $a$ is $\perp$-isotropic if and only if $a \perp A.$
\end{lemma}
\begin{proof} $a \perp A_1'(a),$ by definition.
  Thus the assertion $a\perp  A_1(a)$ is equivalent to  $a \perp A.$
\end{proof}
\begin{remark}\label{orth2}
Orthogonality is a promising tool, with $A^\perp$ contained in every Frobenius radical. On the other hand, one could conceivably have examples for which $A^\perp =A.$
\end{remark}

     \begin{Def}
        An axial algebra $A$ is $X$-anisotropic if every axis in $X$ is $\perp$-anisotropic.
     \end{Def}

 \subsection{Ideals of $A$}

\begin{lemma}(as in \cite[Lemma 5.1]{MSZ})\label{perp2}
    Suppose $\mathcal I \trianglelefteq A$ and let $A^\circ$ be the axial algebra generated by all axes $a$ such that $  \mathcal I\subseteq A_1'(a).$ (For $a$ primitive this means $a\notin \mcI.$) Then
    \begin{enumerate}\eroman
        \item $A^\circ \perp \mathcal I$.
        \item $\mathcal I A^\circ \perp A.$
    \end{enumerate}
\end{lemma}
\begin{proof} (i)
    First we prove by induction that $ a_1\dots a_t \perp I$ for any $t\ge 1.$ If $t = 1$ this is true since $a_1 \in A_1(a_1)$. For $t>1,$ we have $a_t \mathcal I  \subseteq \mathcal I $ so $a_1 \dots a_{t-1}\perp a_t \mathcal I $ by induction, implying $a_1 \dots a_t \perp  \mathcal I $ by Definition~\ref{orth}(4).
    Continuing the argument shows  $z_1\dots z_m\perp \mcI $ for all monomials $z_1,\dots, z_m.$

    (ii)  $\mathcal I \perp A^\circ A$ implies for all $a_i, $ $\mcI \perp (a_t\dots a_1)A$, so $\mcI a_1\dots a_t\perp A$. Continuing inductively, for every product $z_1\dots z_m$ of monomials, $\mcI z_1\dots z_m\perp A$;  hence $\mathcal I A^\circ \perp A.$ \end{proof}

 \subsubsection{The core of a submodule}

\begin{remark}\label{core1}
    Suppose $\mcN$ is a submodule of $\mcA.$ Then the sum of any two submodules of  $\mcN$ is contained in $\mcN$. Hence the sum of all ideals of $A$ contained in $\mcN$ is the unique largest ideal of $A$  contained in $\mcN,$ which we denote as $\core(\mcN)$.
\end{remark}

\begin{lemma}
    $\core(A_1'(a))$ is the unique largest ideal of $A$   intersecting $A_1(a)$ trivially.
\end{lemma}
\begin{proof}
  Obviously $A_1'(a)$, being $a$-homogeneous,  is the largest $C[\gl]$-submodule of $A$ intersecting $A_1(a)$ trivially,  so we take the core.
\end{proof}

\begin{lemma}
    For any submodules $\mcN_i$ of $A,$ $\core(\cap \mcN_i)= \bigcap \core(\mcN_i).$
\end{lemma}
\begin{proof}
    Clearly $\bigcap \core(\mcN_i)$ is an ideal contained in $\cap \mcN_i$, so is contained in $\core(\cap \mcN_i).$ For the other direction, $\core(\cap \mcN_i)$ is contained in each $\core(\mcN_i),$ and thus in $\bigcap \core(\mcN_i).$
\end{proof}

\subsubsection{Blocks}

Following \cite{MSZ}, we define the \textbf{block} $\mathcal I _a$ to be the ideal of $A$ generated by $a,$ and say  $a$ \textbf{dominates} an axis $a'$ if $\mathcal I _{a}\supseteq \mathcal I _{a'}.$

\begin{lemma}\label{comp2}
    $A_0'(a) \subseteq \mathcal I _a,$ for each axis $a.$
\end{lemma}
\begin{proof} $\mcI_a,$ being an ideal, is $a$-homogeneous.
    For any $y= \sum y_\gl \in A,$ $\sum \gl y_\gl =ay\in \mathcal I _a,$ implying $y_\gl \in \mathcal I _a,$ by homogeneity. Thus $A_0'(a) \subseteq \mathcal I _a.$
\end{proof}

\begin{cor}\label{comp3}  Either $\mathcal I _a \cap A_0(a) \ne 0,$  or $\mathcal I _a =A_0'(a).$
\end{cor}

\begin{lemma}
    For any $a\in X$, and any ideal $\mathcal I $ of $A,$   $\pi_a$ restricts to a projection from~$\mathcal I $ to $\mathcal I \cap \mcA_0(a) $.
\end{lemma}
\begin{proof} For $y\in \mathcal I ,$
    $ay = \sum _{0\ne\gr \in S} \gr y_\gr$, implying $y_\gr \in \mathcal I $ for each $\gr\ne 0$, so $\pi_a(y) =y_0\in \mathcal I  \cap \mcA_0(a).$
\end{proof}

\begin{lemma}\label{eq1} The following are equivalent for $a,a' \in X:$\begin{enumerate}\eroman
    \item
    $a$ dominates $a'$
    \item $\pi_a(a')\in \pi_a(\mcI_{a}) .$
    \item  $\pi_a(\mcI_{a'})\subseteq \pi_a(\mcI_{a})$.

\end{enumerate}
\end{lemma}
\begin{proof} $((i)\Rightarrow(ii))$  Obvious.

    $((ii) \Rightarrow(i))$ By assumption, $a'_0 \in \pi_a(\mcI_a),$ so $\sum_{\gr\in S} a'_\gr \in \pi_a(\mcI_{a}) +A_0'(a) \subseteq \mathcal I _a.$

    $((ii) \Leftrightarrow(iii))$ Immediate.
\end{proof}

 \begin{lemma}\label{eq2} When $C$ is a field,
     $a'$ dominates a primitive axis $a$ iff $\mcI_{a'} \not \subseteq A_1'(a).$
 \end{lemma}
\begin{proof}
    $(\Leftrightarrow)$ Each condition says $y_1\ne 0$ for some $y\in \mcI_{a'}.$
\end{proof}

\subsubsection{Digression: Inner ideals}

\begin{Def}
    For   $y_1,y_2 \in A$, define the operator $U_{y_1,y_2}:A\to A$ by $$U_{y_1,y_2}(z) = y_1(y_2z) + y_2(y_1z) - (y_1y_2)z .$$

    An \textbf{quadratic ideal} of $A$ is a subspace $Q$ such that $U_{y_1,y_2}(A)\subseteq Q$ for every $y_1,y_2\in Q$.
 An \textbf{inner ideal} of $A$ is a quadratic ideal $Q$ such that $y^2 \in Q$ for all $y\in Q.$
 \end{Def}

    \begin{example}
        By Remark~\ref{Ua},  for any  axis $a$ in an axial algebra $A$ of Jordan type $\half$, $Fa$ is a quadratic ideal of $A$.
    \end{example}

    \begin{remark}
        This definition specializes to the definition of a quadratic ideal of a linear Jordan algebra, seen by linearizing the quadratic operator, noting that $U_{a,a}= 2U_a,$ being its polarization. Jacobson~\cite[Chapter~4]{J2} has a theory of linear Jordan algebras satisfying the minimum condition on quadratic ideals, so one is led to try to generalize this to axial algebras of Jordan type. So far we have been unable to get useful results along these lines.
    \end{remark}

 \section{Radicals}

In this section, we describe various radicals of an axial algebra $A.$  The definition often depends on the choice of the generating set $X.$ From now on we assume that $C$ is a field $F,$ to simplify the discussion.

\subsection{$A^\perp$}

We already defined $A^\perp$ (with respect to a choice of $X$), which is contained in the Frobenius radical. Hence $A^\perp=0$ when a nonsingular Frobenius form exists.

\begin{lemma}
    $(A/A^\perp)^\perp = 0.$ In particular $A/A^\perp$ is $X$-anisotropic, and $A^\perp$ contains all isotropic axes of $X$.
\end{lemma}
\begin{proof}
    Follows at once from Lemma~\ref{perp0}.
\end{proof}

\subsection{The axial radical}


\begin{Def}\label{MSZ1} Suppose $A$ is an axial algebra, with generating set $X$.
  \begin{enumerate}
      \item   The \textbf{axial radical} $R(A;X)$ is the largest ideal contained in $ A_1'(a)$ for each $a\in X$.
  \end{enumerate}
\end{Def}

\begin{remark}\label{MSZ2}
         $R(A;X)$ is defined, being the core of $\cap _{a\in X} A_1'(a),$ and is the most robust definition of an axial radical. However, it depends on the choice of $X.$ Clearly $R(A;X)\supseteq R(A;X')$ when $X\subseteq X'$, so we can define the ``full'' radical $R(A) =   R(A;\hat X)$, taken over the set $\hat X$ of all $S$-axes.  Thus $R(A) = \cap R(A;X)$ taken over all generating sets $X.$ Likewise one could define  $R_{\operatorname{prim}}(A)  = \cap R(A;X)$ taken over all  generating sets $X$ of primitive $S$-axes.
\end{remark}

The definition of $R(A;X)$ corresponds to \cite{KMS,MSZ} when $A$ is primitive axial over a field, since $y\in \mcI\setminus A_1'(a)$ means $y_1\ne 0$, and thus $a = \varphi_a(y)^{-1}y_1 \in \mcI.$

\begin{lemma}
    $R(A;X)$ is the intersection of the cores of maximal  $X$-axially modular submodules of $\mcA.$
\end{lemma}
\begin{proof}
    As in \cite[Proposition~4.2]{MSZ}.
\end{proof}

\subsubsection{Skeletal elements}

The following notion may serve as a substitute for quasi-invertibility.

\begin{Def} $Z(A;X) = \cap_{a\in X} A_1'(a),$ the intersection of the maximal modular submodules of $A.$
    An element $y\in A$ is \textbf{skeletal} if $y\in Z(A;X),$ i.e. the component $y_1$ with respect to each axis $a\in X$ is $0$.

    An ideal $\mcI$ is \textbf{skeletal} if $\mcI\subseteq Z(A;X)$.
\end{Def}

\begin{lemma}
    $R(A;X)$ contains every skeletal ideal $\mcI$.
\end{lemma}
\begin{proof}
    $\mcI \subseteq Z(A;X)$, so $\mcI = \core(\mcI)\subseteq \core(Z(A;X)).$
\end{proof}

\begin{lemma}
    \label{quot} $R(\bar A;\bar X) = 0$, where $\bar A = A/R(A;X)$.
\end{lemma}
\begin{proof} Write $R(\bar A;\bar X) = \mcI/R(A;X),$ for $R(A;X)\subseteq \mcI \trianglelefteq A.$
    If   $y\in \mcI,$ then, with respect to each $\bar a \in \bar X$, $ 0 = (\bar y)_1 = \overline{y_1},$ implying $y_1 =0$ since  $R(A;X)$ is skeletal. Hence $\mcI = \core(\mcI)\subseteq \core(Z(A;X))= R(A;X).$
\end{proof}

\subsubsection{Applying the orthogonality relation}

\begin{lemma}\label{perp4} For any axial algebra $A$ with generating set $X,$
\begin{enumerate}\eroman
    \item $R(A;X) \subseteq  A^\perp  $.

    \item $R(A;X)= A^\perp$ if $A$ is $X$-anisotropic.
\end{enumerate}

     \end{lemma}
     \begin{proof} (i)  $a\perp Z(A;X) $ for every $a\in X$ since $Z(A;X)\subseteq A_1'(a)$, and hence $z\perp R(A;X)$ for every monomial $z$. Hence $  R(A;X) \subseteq A^\perp$ by Lemma~\ref{perp3}(i).

      (ii) If $R(A;X)\ne A^\perp$ then the ideal $A^\perp$ would contain an axis $a\in X$, contrary to Lemma~\ref{sing}. Hence $R(A;X) = A^\perp $.
     \end{proof}

     Thus, the  $X$-anisotropic condition is very powerful, and enables us to apply  the theory of \cite{MSZ}.

\begin{prop}(as in \cite[Proposition~5.3]{MSZ})
    Suppose $A$ has an orthogonality relation, with respect to a set $X$ of primitive axes, with $A^\perp = 0,$ and take any ideal $\mathcal I $ of $A.$ Let $X_1 = X \cap \mathcal I $ and $X_2 = X \setminus X_1.$ Then $\mathcal I $ is the axial algebra $A_1$ generated by $X_1$,  and is perpendicular to the axial algebra $A_2$ generated by $X_2$. In other words $A = A_1 \oplus A_2,$ so $\mcI \cap A_2$ contains no axis of $X$ (since $X_2 \cap \mcI = \emptyset$). Hence $\mcI \cap A_2 \subseteq R(A;X)=0,$ implying $\mcI \subseteq A_1 \subseteq \mcI,$ implying $\mcI = A_1.$
\end{prop}
\begin{proof}   $R(A;X) = 0$ by Lemma~\ref{perp4}, so every nonzero ideal contains an axis of $X$. Also
  $y\perp z$  for $y\in \mathcal I _1$ and $z\in X_2$, by Lemma~\ref{perp2}, which also  implies $AA_1 \perp A_2 ,$ implying  $A_1A_2\subseteq A^\perp= 0.$ But $X\subseteq A_1\cup A_2,$ implying $A = A_1+ A_2.$ Hence both $A_1,A_2\trianglelefteq A$.
  $(A_1 \cap A_2)^2 =0.$ Hence $A_1 \cap A_2$ does not contain any nontrivial idempotent, so is contained in $R(A;X)=0.$
\end{proof}

We are now in a position to apply  \cite{MSZ}. We say that two axes $a$ and $a'$ are \textbf{equivalent} if the blocks $\mcI_{a} = \mcI_{a'}$. Then by the proof of \cite[Lemma~6.1 and Proposition 6.2]{MSZ}, mutatis mutandis,  we could show that if $A^\perp =0$, then $A = \oplus_{a\in X'}I_a,$
taken over a transversal $X'$ of axes in $X$.

 \subsection{The $\gr$-Baer radical}

  Lemma~\ref{normal} permits an alternate version of the radical of a $\gr$-axial algebra, without using a Frobenius form.

Since $Z(A;X)$ can be described directly,  a Baer-type construction  may
be useful to equalize $Z(A;X)$ with $R(A;X)$.

   \begin{Def}\label{gen7}
  The $\gr$-\textbf{Baer radical} $\mathcal{R}(A;X)$ is defined inductively as follows:
$\mathcal{R}_0$ is the ideal generated  by $ Z(A;X),$ and $Z_0 = Z(A;X).$

For any  $i\in \N,$ letting $\bar A = A/\mathcal{R}_i,$  $\mathcal{R}_{i+1}\supseteq \mathcal{R}_i$ with    $\mathcal{R}_{i+1}/\mathcal{R}_i$ is the ideal generated by $ Z(\bar A;\bar X) $, where $\bar X$ is the set of nonzero images of axes of $X$. $\mathcal{R}(A;X)= \cup_{i\ge 0} \mathcal{R}_i.$
\end{Def}
Let us see that for $X$ finite, this process must terminate.

\begin{lemma} Notation as above, for $\bar A = A/\mathcal{R}_i(A)$,
\begin{enumerate}\eroman
    \item
    If $ Z(\bar A;\bar X)$ is already an ideal then $  \mathcal{R}_{i+1}=\mathcal{R}_i,$ so  the process stops.
    \item If $ Z(\bar A;\bar X)$ is not an ideal, then the ideal $\mathcal{R}_{i+1}/\mathcal{R}_{i}$ that it generates properly  contains $\core( Z(\bar A;\bar X)) = R(\bar A, \bar X)$, and thus contains an axis of $\bar X,$ which is absorbed in passing to $A/\mathcal{R}_{i+1}.$
  \item  The number of stages in the inductive process does not exceed $|X|.$
\end{enumerate}
\end{lemma}
\begin{proof}
    (i) Clearly $ Z(\bar A;\bar X) =R(\bar A;\bar X),$ so $\mathcal{R}_{i+1}=\mathcal{R}_{i}.$

    (ii) By definition, $R(\bar A;X) = \core(\bar Z) \ne \bar Z,$ so $\mathcal{R}_{i+1}$ contains a nonzero axis in $\bar X,$ whose image becomes 0.

    (iii) Apply (ii) inductively.
\end{proof}

 By construction, $\mathcal{R}(A;X) \supseteq Z(A;X) \supseteq R(A;X),$ equality holding for $C$ a field and $X$ primitive if and only if $X\cap \mathcal{R}(A;X) =\emptyset $.

 \subsection{Special case: The radical of a Jordan algebra}

\begin{Def}{[following \cite[Exercise~1.11.5, p.~55]{J2}]} An element $y\in A$ is \textbf{quasi-invertible}\footnote{McCrimmon~\cite{McC} more generally defined quasi-invertibility for quadratic Jordan algebras, but we stay with the linear case since we are assuming that $\half\in F.$} with \textbf{quasi-inverse} $z$ if $y^2z - y^2 =yz = y+z $.  An ideal is \textbf{quasi-invertible} if it consists of  quasi-invertible elements.
\end{Def}

Clearly a quasi-invertible ideal cannot contain an idempotent, so is contained in $R(A;X)$ when $X$ consists of primitive axes. On the other hand, $R(A)$ could conceivably contain idempotents, so is not quasi-invertible in general.

Jacobson~\cite[Exercise 1.11.6]{J2} asserts that  in a Jordan algebra  all of whose elements are algebraic, an ideal is quasi-invertible if and only if it is nil. This holds for Clifford Jordan algebras and for locally finite dimensional Jordan algebras. McCrimmon~\cite{McC}  proved that the sum of quasi-invertible ideals is quasi-invertible, and thus yields a maximal quasi-invertible ideal, which we call  the \textbf{quasi-invertible radical} $QI(A).$

Hogben and McCrimmon \cite{HM} defined $a$-modular inner ideals with respect to a modulus $a$, calling a modular inner ideal maximal if it is maximal $a$-modular for some~$a.$ They then showed that $QI(A)$ is the intersection of the maximal modular inner ideals, and thus the intersection of the maximal modular ideals.

\subsubsection{Isotopes}

This subsection is largely negative. One of the main tools in Jordan theory is the $a$-\textbf{homotope} $(A,a)$ of a Jordan algebra, given by multiplication $a_1 \cdot_a a_2 = \{ a_1 a a_2\}.$ For a primitive axis $a,$ $a\cdot_b a = 2a(ab)-ab = \sum_{\gl \in S} (2\gl^2 -\gl)b_\gl,$ which is $a$ precisely when $2\gl^2 -\gl =0$ (i.e., Jordan type $ \{\half\}$) and $\varphi_a(b)=1,$ so this seems to be useful only for Jordan type $\half.$
\section{Some further questions}
\begin{ques}
 \begin{enumerate}
     \item    Under what conditions is $Z(A;X)$ an ideal?

       \item    How can we determine if $A$ is $X$-anisotropic, in the absence of a Frobenius form?  In particular,  are axial algebras of Monster type $X$-anisotropic, for a suitable generating set $X$?

             \item Does there exist a nontrivial example for which $A^\perp = A?$

            \item  Under what conditions does $R(A;X)=R(A)$?

             \item Does $R(A;X)$ contain every  skeletal quadratic ideal of $A$?

               \item If the answer to (5) is negative, what is the ``correct'' analog of inner ideal in axial algebras (perhaps of Jordan type, so that one could use Miyamoto involutions)?
 \end{enumerate}
\end{ques}

One promising class for closer scrutiny may be the \textbf{quasi-definite primitive axial algebras of Jordan type $\half$} introduced by \cite{GGu}, defined by $\varphi_a(b)\ne 1$ for any axes $a\ne b.$

\end{document}